\documentclass[reqno]{amsart}
\usepackage{amssymb}
\usepackage{amsmath}
\usepackage{amsfonts}
\usepackage{hyperref}

\theoremstyle{plain}

\newtheorem{corollary}{Corollary}

\newtheorem{definition}{Definition}
\newtheorem{example}{Example}

\newtheorem{lemma}{Lemma}

\newtheorem{theorem}{Theorem}
\numberwithin{equation}{section}

\begin{document}
\title[Linear combination of harmonic functions]{quantum approach on linear combination of harmonic univalent mappings 
}
\author[ O. Mishra and S. Porwal]{
	 Omendra Mishra${}^{1}$ and Saurabh Porwal${}^{2}$}
\date{}
\maketitle
\begin{center}
	${}^{1}$Department of Mathematics\\ Rajat P.G. College \\Lucknow-226007, (U.P.), India.\\
	e-mail: mishraomendra@gmail.com; https://orcid.org/0000-0001-9614-8656.\\
${}^{2}$Department of Mathematics\\
Ram Sahai Government Degree College \\Bairi-Shivrajpur, Kanpur-209205, (U.P.), India.\\
e-mail: saurabhjcb@rediffmail.com; https://orcid.org/0000-0003-0847-3550.\\
\end{center}
\vspace{.2in}
\maketitle
\begin{abstract}
In this paper using $q$ calculus operator we obtain  some sufficient conditions on $f_1$ and $f_2$
 so that their linear combination $%
 f=tf_{1}+(1-t)f_{2},\ t\in \left[ 0,1\right] $, is univalent and convex in
 the direction of the real axis. Some examples are also illustrated to support our main results.
\end{abstract}


\maketitle



\section{Introduction}

The theory of $q$-calculus has motivated the researchers due to its
applications in the field of physical sciences, specially in quantum
physics. Jackson \cite{jack2,FJ} was the first to give some applications of $%
q$-calculus by introducing the $q$-analogues of derivative and integral.
Jackson's $q$-derivative operator $\partial _{q}$ on a function $h$ analytic
in $\mathbb{D}=\left\{ z\in \mathbb{C}:\left\vert {z}\right\vert <1\right\} $
is defined for $0<q<1,$ by 
\begin{equation*}
\partial _{q}h(z)=\left\{ 
\begin{array}{c}
\frac{h(z)-h(qz)}{(1-q)z}\ \ z\neq 0, \\ 
h^{\prime }(0)\quad z=0.%
\end{array}%
\right. \ 
\end{equation*}%
 
For a power function $h(z)=z^{k},\quad k\in \mathbb{N}=\{1, \: 2, \: 3, \: \cdots\},$ 
\begin{equation*}
\partial _{q}h(z)=\partial _{q}(z^{k})=[k]_{q}z^{k-1},
\end{equation*}%
 
$\int_{0}^{z}t^nd_qt=z(1-q)\sum_{k=0}^{\infty}q^k(zq^k)^n=\frac{z^{n+1}}{[n+1]_q}$.\\
where $[k]_{q}$ is the $q$-integer number $k$ defined by 
\begin{equation*}
\lbrack k]_{q}=\frac{1-q^{k}}{1-q}=1+q+q^{2}+...q^{k-1}, 
\end{equation*}%
For any non-negative integer $k$ the $q$-number factorial is defined by
\begin{equation*}
	[k]_{q}!=[1]_{q}[2]_{q}[3]_{q}\dots[k]_{q}\qquad([0]_{q}!=1).
\end{equation*} 
For more detailed study see \cite{V}. Clearly, $\underset{q\rightarrow 1^{-}}%
{\lim }[k]_{q}=k$ and $\underset{q\rightarrow 1^{-}}{\lim }\partial
_{q}h(z)=h^{\prime }(z).$
Research work in connection with function theory and $q$-calculus was first
introduced by Ismail \textit{et al.} \cite{Is}. Recently, $q$-calculus is
involved in the theory of analytic functions in the work \cite%
{GoSiva,HussKhanDarus,mohmSokol}. But research on $q$-calculus
in connection with harmonic functions is fairly new and not much published
(one may find papers \cite{jm}, \cite{MuruJahan}, \cite{spag}, \cite{Bravinder}, \cite{Bravinder 2},  \cite{pso}, \cite{srietal}).

 Let $\mathcal{H}$ denotes the class of complex-valued functions $%
 f=u+iv$ which are harmonic in the unit disk $\mathbb{D}=\left\{ z\in \mathbb{%
 	C}:\left\vert {z}\right\vert <1\right\} ,$ where $u$ and $v$ are real-valued
 harmonic functions in $\mathbb{D}$. Functions $f$ ${\in \mathcal{H}}$ can
 also be expressed as $f=h+\overline{g},$ where $h$ and $g$ are analytic in $%
 \mathbb{D},$ called the analytic and co-analytic parts of $f$, respectively.
 The Jacobian of $f=h+\overline{g}$ is given by $J_{f}(z)=|{h^{\prime }(z)}%
 |^{2}-|{\ g^{\prime }(z)}|^{2}$.
 
 According to the Lewy's Theorem, every harmonic function $f=h+\overline{g}%
 \in \mathcal{H}$ is locally univalent and sense preserving in $\mathbb{D}$
 if and only if $J_{f}(z)>0$ in $\mathbb{D}$ which is equivalent to the
 existence of an analytic function $\omega (z)=g^{\prime }(z)/h^{\prime }(z)$
 in $\mathbb{D}$ such that
 \begin{equation*}
 	|{\omega }(z)|<1\ \ \ \text{for all \ }z\in \mathbb{D}.
 \end{equation*}%
 The function $\omega$ is called the dilatation of $f$. By requiring harmonic
 functions to be sense-preserving, we retain some basic properties exhibited
 by analytic functions, such as the open mapping property, the argument
 principal, and zeros being isolated (see for detail \cite{PWR}). \ The class
 of all univalent, sense preserving harmonic functions $f=h+ \overline{g} \in
 \mathcal{H}$, with the normalized conditions $h(0)=0=g(0)$ and $h^{\prime
 }(0)=1$ is denoted by $S_{\mathcal{H}}$. If the function $f=h+\overline{g}%
 \in S_{\mathcal{H}},$ then $h$ and $g$ are of the form
 \begin{equation*}
 	h(z)=z+\sum_{n=2}^{\infty }a_{n}z^{n}\text{ \ \ and \ \ \ }
 	g(z)=\sum_{n=1}^{\infty }b_{n}z^{n}\text{ \ \ }\left( \left\vert
 	b_{1}\right\vert <1;z\in \mathbb{D}\right) .  
 \end{equation*}%
 A subclass of functions $f=h+\overline{g}\in S_{\mathcal{H}}$ with the
 condition $g^{\prime }(0)=0\ \left( \text{or }\omega _{f}(0)=0\right) $ is
 denoted by $S_{\mathcal{H}}^{0}.$\ Further, the subclasses of functions $f$
 in $S_{\mathcal{H}}$ $\left( S_{\mathcal{H}}^{0}\right)$, denoted by $K_{
 	\mathcal{H}}$ $\left( K_{\mathcal{H}}^{0}\right), $ consists of functions $f$
 that map the unit disk $\mathbb{D}$ onto a convex region.

\begin{definition} \cite{MR 12}	\textbf{Harmonic right half-plane mappings}
	\newline
	A mapping $f=h+\bar{g}$ is said to be a harmonic right half-plane mapping if it maps $\mathbb{D}$ onto a right-half plane
	\begin{equation*}
	H_0=\{w:\Re e(w)>-1/2\}.
	\end{equation*}
	For
	\begin{equation}
	h(z)+g(z)=\frac{z}{1-z} \qquad \text{and} \quad \omega(z)=\frac{g'(z)}{h'(z)}=-z\nonumber
	\end{equation} 
	\begin{equation}
	h'(z)+g'(z)=\frac{1}{(1-z)^2} \qquad g'(z)=-zh'(z)\nonumber
	\end{equation} 
	we get
	\begin{equation}
	h'(z)=\frac{1}{(1-z)^3} \qquad g'(z)=\frac{-z}{(1-z)^3}\nonumber
	\end{equation}
	which on  integration and normalization gives 
	\begin{equation}
	h(z)=\frac{z-\frac{z^2}{2}}{(1-z)^2} \qquad g(z)=\frac{\frac{-z^2}{2}}{(1-z)^2}.\nonumber
	\end{equation}
	Thus the mapping $f=h+\bar{g}$ given by
	\begin{equation}
	f=\frac{z-\frac{z^2}{2}}{(1-z)^2} -\overline{\frac{\frac{z^2}{2}}{(1-z)^2}}.\label{q half}
	\end{equation}
	\\
\end{definition}
The class of
right half-plane mappings $f\in S^{0}(H_{0})$ that map $\mathbb{D}$ onto $f(%
\mathbb{D})=H_{0}=\{w:\Re e(w)>-1/2\}$, and such a mapping clearly assumes
the form
\begin{equation}
	h(z)+g(z)=\frac{z}{1-z}.  \notag
\end{equation}%

We now recall fundamental result, called shearing theorem, due to Clunie and
Sheil-Small as follows:

\begin{theorem}
	\label{thA}\cite{JT 84}A locally univalent harmonic function $f=h+\overline{g}$
	in $\mathbb{D}$ is a univalent mapping of real axis $\mathbb{D}$ onto a domain convex
	in the direction of real axis if and only if $h-g$ is a
	analytic univalent mapping of $\mathbb{D}$ onto a domain convex in the
	direction of real axis.
\end{theorem}
In this paper using quantum approach we  find the new harmonic functions.
Such that, if we have an analytic function of the form
\begin{example}
		$h(z)-g(z)=z-\frac{1}{2}z^2$ \qquad \text{and} \quad $\omega_q(z)=\frac{\partial _{q}(g(z))}{\partial _{q}(h(z))}=\frac{[2]_q}{2}z$ \nonumber 
	\begin{equation}
	\partial _{q}(h(z))-\partial _{q}(g(z))=1-\frac{[2]_q}{2}z \qquad \partial _{q}(g(z))=\frac{[2]_q}{2}z\partial _{q}(h(z))\nonumber
	\end{equation} 
	we get
	\begin{equation}
	\partial _{q}(h(z))=1\nonumber
	\end{equation}
	which on  q-integration and normalization gives 
	\begin{equation}
	h(z)=z \qquad g(z)=\frac{[2]_q}{4}z^2\nonumber
	\end{equation}
	Thus the mapping $f=h+\bar{g}$ given by
	\begin{equation}
	f=z-\frac{[2]_q}{4}\overline{z^2}\nonumber
	\end{equation}
	On $q\rightarrow 1^{-}$, we get the original function $	f=z-\frac{1}{2}\overline{z^2}$
\end{example} 
	\begin{definition}
		\textbf{Harmonic q right half-plane mappings}: \\
		If we have an analytic function of the form
		\begin{equation}
		h(z)+g(z)=\frac{z}{1-z} \qquad \text{and}\qquad \omega_q(z)=\frac{z^2-2z+qz}{1-qz}.  \notag
		\end{equation}%
			Then
    		\begin{equation*}
    		\partial _{q}(h(z))+\partial _{q}(g(z))=\frac{1}{(1-z)(1-qz)}
    		\end{equation*}
    			we get
    			\begin{equation}
    			\partial _{q}(h(z))=\frac{1}{(1-z)^3}\nonumber
    			\end{equation}
    			which on  q-integration and normalization gives 
		\begin{equation*}
		h(z)=\sum_{n=0}^{\infty }\frac{(n+1)(n+2)}{2[n+1]_q}z^{n+1} 
		\end{equation*}
		Using this we get
		\begin{equation*}
		g(z)=\sum_{n=0}^{\infty }\frac{2[n+1]_q-(n+1)(n+2)}{2[n+1]_q}z^{n+1} 
    		\end{equation*}
		Thus $h$ and $g$ of the form
		\begin{equation*}
			h(z)=\sum_{n=0}^{\infty }\frac{(n+1)(n+2)}{2[n+1]_q}z^{n+1} \qquad 	g(z)=\sum_{n=0}^{\infty }\frac{2[n+1]_q-(n+1)(n+2)}{2[n+1]_q}z^{n+1} 
		\end{equation*}
		On $q\rightarrow 1^{-}$, we get the mapping defined in \eqref{q half}.
			\end{definition}
	
		\begin{lemma}
			\label{lemma2}\cite{JT 84} Let $\Omega \subset\mathbb{C}$ be a domain convex
			in the direction of the real axis. Also let $p$ be a real-valued continuous
			function in $\Omega .$ Then the mapping $\omega \mapsto \omega + p(\omega )$
			is univalent in $\Omega $ if and only if it is locally univalent. If it is
			univalent, then its range is convex in the direction of the real axis.
		\end{lemma}
		\begin{lemma}
			\label{lemma1}\cite{Pm}Let $f$ be analytic function in $\mathbb{D}$ with $%
			f(0)=0$ and $f^{\prime }(0)\neq 0.$ Suppose also that%
			\begin{equation*}
				\varphi (z)=\frac{z}{\left( 1+ze^{i\theta }\right) \left( 1+ze^{-i\theta
					}\right) }\ \ \left( \theta \in\mathbb{R};z\in \mathbb{D}\right) .
			\end{equation*}%
			If%
			\begin{equation*}
				\Re \left( \frac{zf^{\prime }(z)}{\varphi (z)}\right) >0\ \ \left( z\in 
				\mathbb{D}\right) ,  
			\end{equation*}%
			then $f$ is convex in the direction of real axis.
		\end{lemma}
Dorff and Rolf \cite{MR 12} applied another way of constructing a univalent
harmonic map by taking two suitable harmonic maps $f_{1}$ and $f_{2}$ with
same dilatations, whose linear combination $f_3=tf_{1}+(1-t)f_{2},\ t\in \left[
0,1\right] $ is univalent and convex in the direction of the imaginary axis.
Wang et al. \cite{wang} derived several sufficient conditions on harmonic
univalent functions $f_{1}$ and $f_{2}$ so that their linear combination $%
f=tf_{1}+(1-t)f_{2},\ t\in \left[ 0,1\right] $, is univalent and convex in
the direction of the real axis. More results on the linear combination $f$
of $f_{1}$ and $f_{2}$ may also be found in \cite{salas,kumar
	gupta, psomen,shiwang,sun rasila,wang} etc. (also see the references cited in these). In this paper using quantum approach we find sufficient conditions on $f_1$ and $f_2$
so that their linear combination $%
f=tf_{1}+(1-t)f_{2},\ t\in \left[ 0,1\right] $, is univalent and convex in
the direction of the real axis. Some examples are also illustrated to support our main results.

\section{main results}
First, we give $q$-analogue of Lemma \ref{lemma1}
	\begin{lemma}
		\label{c2lemma1}Let the function $f:\mathbb{D}\rightarrow C$ be
		an analytic function with $f(0)=0$ and $\partial _q(f(0))\neq 0$. Suppose that
		\begin{equation}
		\underset{q\rightarrow 1^{-}}%
		{\lim }\varphi_q (z) =   	\varphi (z)=\frac{z}{\left( 1+ze^{i\theta }\right) \left( 1+ze^{-i\theta
			}\right) }\ \ \left( \theta \in \mathbb{R};z\in \mathbb{D}\right) .
		\label{c2varphi}
		\end{equation}%
		If the function $f$\ satisfy
		\begin{equation}
		\Re e\left( \frac{z\partial_q(f(z))}{\varphi_q     (z)}\right) >0\ \ \left( z\in
		\mathbb{D}\right) ,  \label{c2conlemma1}
		\end{equation}%
		then the function $f$ is convex in the direction of real axis.
	\end{lemma}
\begin{theorem}
	\label{th1}Let for $j=1,2,f_{j}=h_{j}+\overline{g_{j}}\in S_{\mathcal{H}}$
	with%
	\begin{equation}
	F_j=	h_{j}(z)-g_{j}(z) .  \label{phij}
	\end{equation}%
	  If ${\omega }_{q_{1}}{%
		(z)=\omega }_{q_{2}}{(z)}$ and  satisfy the condition $\Re \left( 
	\frac{z\partial_q(F_j (z))}{\varphi_q(z)}\right) >0$ for some function $%
	\varphi_q(z)$ given by \eqref{c2conlemma1}, then $f_3=tf_{1}+(1-t)f_{2},\ t\in \left[ 0,1\right]$ is univalent and  convex
	in the direction of real axis.
    \end{theorem}

\begin{proof}
 Since, for $j=1,2,{\omega }%
	_{q_{j}}{(z)}=\frac{\partial _q(g_{j}(z))}{\partial _q(h_{j}(z))},$ with $%
	\left\vert {\omega }_{q_{j}}{(z)}\right\vert <1.$ If ${\omega }_{q_{1}}{(z)=\omega }_{q _{2}}{(z)}$ then
	\begin{align*}
		 {\omega }_{q_3}{(z)}&= \frac{t\partial _q(g_{1}(z))%
			+(1-t)\partial _q(g_{2}(z))}{t\partial _q(h_{1}(z))+(1-t)\partial _q(h_{2 }(z))}%
		  \notag \\
		&= \frac{t{\omega }%
			_{q_{1}}{(z)}\partial _q(h_{1}(z))%
			+(1-t){\omega }%
			_{q_{1}}{(z)}\partial _q(h_{2}(z))}{t\partial _q(h_{1}(z))+(1-t)\partial _q(h_{2 }(z))} \\
	&=\left\vert {\omega }_{q_{1}}{(z)}\right\vert <1. \notag
\end{align*}%
 To show $f_3$ is convex in the direction of real axis, we have to show by Lemma \eqref{c2lemma1}.
\begin{align*}
\Re \left( \frac{z(\partial _q(h_{j}(z))-\partial _q(g_{j}(z)))}{\varphi_q(z)}\right)>0
  \end{align*}
 Now
 \begin{align*}
& \Re \left( \frac{z(\partial _q(h_{3}(z))-\partial _q(g_{3}(z)))}{\varphi_q(z)}\right)&
 \\&=\Re\left(\frac{z}{\varphi_q(z)}\left[t(\partial _q(h_{1}(z))-\partial _q(g_{1}(z))+(1-t)(\partial _q(h_{2}(z))-\partial _q(g_{2}(z))\right]\right)\\
 &=t\Re\left(\frac{z}{\varphi_q(z)}(\partial _q(h_{1}(z))-\partial _q(g_{1}(z))\right)+(1-t)\Re\left(\frac{z}{\varphi_q(z)}(\partial _q(h_{2}(z))-\partial _q(g_{2}(z))\right) \\
 &>0.
  \end{align*}
 \end{proof}
A generalization of Theorem \ref{th1} may be given as follows:
\begin{corollary}
	\label{co1}Let for $j=1,2,...,n,f_{j}=h_{j}+\overline{g_{j}}\in S_{\mathcal{H%
		}}$ with the condition (\ref{phij}).  If ${\omega }%
		_{q_{1}}{(z)=\omega }_{q_{2}}{(z)=...=\omega }_{q_{n}}{(z)}$ and 
		satisfy the condition $\Re \left( \frac{z\partial_q(F_J (z))}{\varphi_q(z)}%
		\right) >0$ for some function $\varphi_q(z)$ given by \eqref{c2conlemma1}. Then $F=\sum_{j=1}^{n}t_{j}f_{j},\
		t_{j}\in \left[ 0,1\right] $ such that $\sum_{j=1}^{n}t_{j}=1.$ is univalent and convex in the direction of real axis. 
	\end{corollary}
\begin{theorem}
\label{q th}Let for $j=1,2,f_{j}=h_{j}+\overline{g_{j}}\in S_{\mathcal{%
			H}}$ be convex in the direction of the real axis. If $\Re \left\{ \left( 1-{\omega }_{q_{1}}%
	\overline{{\omega }_{q_{2}}}\right) \partial _q(h_{1}(z)\overline{\partial _q(h_{2}(z)}%
	\right\} \geq 0$, then   $f_3=tf_{1}+(1-t)f_{2},\
	t\in \left[ 0,1\right]$ is convex in the direction of real axis. .
\end{theorem}
\begin{proof}
Since, for $j=1,2,{\omega }%
_{q_{j}}{(z)}=\frac{\partial _q(g_{j}(z))}{\partial _q(h_{j}(z))},$ with $%
\left\vert {\omega }_{q_{j}}{(z)}\right\vert <1.$ We get 
\begin{align*}
\vert {\omega }_{q_3}{(z)}\vert &= \left \vert  \frac{t\partial _q(g_{1}(z))%
	+(1-t)\partial _q(g_{2}(z))}{t\partial _q(h_{1}(z))+(1-t)\partial _q(h_{2 }(z))} \right \vert \notag \\
&= \frac{\vert t{\omega }%
	_{q_{1}}{(z)}\partial _q(h_{1}(z))%
	+(1-t){\omega }%
	_{q_{2}}{(z)}\partial _q(h_{2}(z))\vert}{\vert t\partial _q(h_{1}(z))+(1-t)\partial _q(h_{2 }(z))\vert } \\
\end{align*}
To show  $\left\vert {\omega }_{q_{3}}\right\vert <1$. We have to prove%
\begin{eqnarray*}
	&&\left\vert t\partial _q(h_{1}(z))+(1-t)\partial _q(h_{2}(z))\right\vert
	^{2}-\left\vert t{\omega }_{q_{1}}\partial _q(h_{1}(z))+(1-t){\omega }%
	_{q_{2}}\partial _q(h_{2}(z))\right\vert ^{2} \\
	&=&\left( t\partial _q(h_{1}(z))+(1-t)\partial _q(h_{2}(z))\right) \left( t\overline{%
		\partial _q(h_{1}(z))}+(1-t)\overline{\partial _q(h_{2}(z))}\right) \\
	&&-\left( t{\omega }_{q_{1}}\partial _q(h_{1}(z))+(1-t){\omega }%
	_{q_{2}}\partial _q(h_{2}(z))\right) \left( t\overline{{\omega }%
		_{q_{1}}\partial _q(h_{1}(z))}+(1-t)\overline{{\omega }_{q_{2}}\partial _q(h_{2}(z))}%
	\right) \\
&=&t^{2}\left( 1-\left\vert {\omega }_{q_{1}}\right\vert ^{2}\right)
\left\vert \partial _q(h_{1}(z))\right\vert ^{2}+\left( 1-t\right) ^{2}\left(
1-\left\vert {\omega }_{q_{2}}\right\vert ^{2}\right) \left\vert
\partial _q(h_{2}(z))\right\vert ^{2} \\
&&+2t(1-t)\Re \left\{ \left( 1-{\omega }_{q_{1}}\overline{{\omega }_{q_{2}}}%
\right) \partial _q(h_{1}(z))\overline{\partial _q(h_{2}(z))}\right\} \\
&>&0
\end{eqnarray*}%
We use Lemma \ref{lemma2} and Theorem \ref{thA} which shows that, $F_{j}:=h_{j}-g_{j}$ is
univalent in $\mathbb{D}$ and $\Omega _{j}=F_{j}\left( \mathbb{D}\right) $
is convex in the direction of the real axis for each $j=1,2$. We may write $%
f_{j}=F_{j}+2\Re \left( g_{j}\right) $ and 
\begin{equation*}
f_{j}\left( F_{j}^{-1}\left( w\right) \right) =w+2\Re \left( g_{j}\left(
F_{j}^{-1}\left( w\right) \right) \right) =w+q_{j}\left( w\right)
\end{equation*}%
is univalent in $\Omega _{j}$ for each $j=1,2,$ where $q_{j}\left( w\right) $
is real valued continuous function. Let $F:=h-g$ and $\Omega =F\left( 
\mathbb{D}\right) $. Then 
\begin{eqnarray*}
	f\left( F^{-1}\left( w\right) \right) &=&tf_{1}\left( F_{1}^{-1}\left(
	w\right) \right) +\left( 1-t\right) f_{2}\left( F_{2}^{-1}\left( w\right)
	\right) \\
	&=&t\left( w+q_{1}\left( w\right) \right) +\left( 1-t\right) \left(
	w+q_{2}\left( w\right) \right) \\
	&=&w+\left( tq_{1}\left( w\right) +\left( 1-t\right) q_{2}\left( w\right)
	\right) \\
	&=&w+q(w)
\end{eqnarray*}%
is univalent in $\Omega $ which by Lemma \ref{lemma2} is convex in the
direction of the real axis.
	\end{proof}
	\section{Examples}
	In this section, we give some examples based on Theorem \ref{th1} and Theorem \ref{q th} to verify our main results.
	\begin{example}
		Let for $j=1,2,f_{j}=h_{j}+\overline{g_{j}}\in S_{\mathcal{H}}$
		with%
		\begin{equation}
		F_j=	h_{j}(z)-g_{j}(z)=\frac{z}{1-z}.  
		\end{equation}%
		If ${\omega }_{q_{1}}{%
			(z)=\omega }_{q_{2}}{(z)}$ and  $f_3=tf_{1}+(1-t)f_{2},\ t\in \left[ 0,1\right]$ is univalent and  convex
		in the direction of real axis. By Theorem \ref{th1} we have 
		\begin{align*}
		\Re \left( \frac{z(\partial _q(h_{3}(z))-\partial _q(g_{3}(z)))}{\varphi_q(z)}\right)>0
		\end{align*}
Since,	$\partial _q(h_{1}(z))-\partial _q(g_{1}(z))=\frac{1}{(1-z)(1-qz)}$, $\partial _q(h_{2}(z))-\partial _q(g_{2}(z))=\frac{1}{(1-z)(1-qz)}$ and $\varphi_q(z)=\frac{z}{(1-z)(1-qz)}$.        
		
	\end{example}
	\begin{example}
	$f_1=z-\frac{[2]_q}{2}\bar{z^2}$\qquad and \qquad $f_1=z+\frac{[3]_q}{3}\bar{z^2}$
	and  $f_3=tf_{1}+(1-t)f_{2},\
	t\in \left[ 0,1\right]$ with $\omega_{q_1}(z)=-\frac{[2]_q}{2}z$ and  $\omega_{q_2}(z)=\frac{[3]_q}{3}z^2$\\
	$\vert {\omega }_{q_3}{(z)} \vert \leq t\frac{[2]_q}{2}\vert z \vert +(1-t)\frac{[3]_q}{3}\vert z^2\vert<t\frac{[2]_q}{2}\vert z \vert +(1-t)\frac{[3]_q}{3}\vert z\vert <t\frac{[2]_q}{2} +(1-t)\frac{[3]_q}{3}\vert<1$ 
	and
	$\Re \left\{ \left( 1-{\omega }_{q_{1}}%
	\overline{{\omega }_{q_{2}}}\right) \partial _q(h_{1}(z)\overline{\partial _q(h_{2}(z)}%
	\right\}=\Re \left(1+\frac{[2]_q[3]_q}{6}\vert z^2\vert \bar{z}\right)>0.$
	\end{example}

\end{document}